\documentclass{amsart}

\usepackage{amssymb}

\usepackage{graphicx}
\usepackage{psfrag}

\newtheorem{theorem}{Theorem}[section]

\newtheorem{proposition}[theorem]{Proposition}
\newtheorem{corollary}[theorem]{Corollary}

\theoremstyle{definition}

\theoremstyle{remark}
\newtheorem{remark}[theorem]{Remark}

\numberwithin{equation}{section}

\newcommand{\TryPackage}[3]{\IfFileExists{#1.sty}{\usepackage{#1}#2}{#3}}
\TryPackage{mathrsfs}{\renewcommand{\mathcal}{\mathscr}}{%
     \TryPackage{eucal}{}{}}

\newcommand{\ga}{\gamma}

\newcommand{\la}{\lambda}

\newcommand{\Si}{\Sigma}

\newcommand{\CC}{{\mathbb C}}

\newcommand{\QQ}{{\mathbb Q}}

\newcommand{\mer}{{\mathcal M}} 
\newcommand{\lng}{{\mathcal L}}  
\newcommand{\SLC}{{SL_2({\mathbb C})}}
\newcommand{\slc}{{{\mathfrak s}{\mathfrak l}_2({\mathbb C})}}

\newcommand{\Ad}{\operatorname{\it Ad}}

\begin{document}

\title{ Splicing and the SL$_2(\CC)$ Casson invariant}


\author{Hans U. Boden}
\address{Department of Mathematics \& Statistics, McMaster University, Hamilton, Ontario, L8S 4K1 Canada}
\curraddr{} \email{boden@mcmaster.ca}
\thanks{The first named author was supported by a grant from the Natural Sciences and Engineering Research Council of Canada.}

\author{Cynthia L. Curtis}
\address{Department of Mathematics \& Statistics, The College of New Jersey, Ewing, NJ, 08628 USA}
\curraddr{} \email{ccurtis@tcnj.edu}
\thanks{}



\dedicatory{}


\begin{abstract}
We establish a formula for the $\SLC$ Casson invariant of
spliced sums of homology spheres along knots. Along the way,
we show that the $\SLC$ Casson invariant vanishes for
spliced sums along knots in $S^3$.
\end{abstract}

\maketitle


\section{Introduction}
In  \cite{KM1}, Kronheimer and Mrowka prove that all nontrivial knots in $S^3$ have Property P.
Their proof is based on strong existence results for
irreducible $SU(2)$ representations of 3-manifolds obtained by Dehn surgery.
It remains an interesting and important problem to determine whether a given
3-manifold admits irreducible $SU(2)$ representations.
For example, for  homology spheres $\Si$,
nontriviality of the Casson invariant  or Floer homology
implies the existence of an irreducible $SU(2)$ representation. Since every irreducible $SU(2)$
representation is also irreducible as an $\SLC$ representation,
one expects stronger results for $\SLC$. For example,
Boyer and Zhang \cite{BZ3}, and independently Dunfield and Garoufalidis \cite{DG},  show that
any nontrivial knot $K$ in $S^3$ has nontrivial $A$-polynomial
by using \cite{KM2} to establish the existence of an arc of irreducible $\SLC$ characters
on $\pi_1(S^3 \setminus K)$.

Given a closed 3-manifold $\Si$, the $\SLC$ Casson invariant
$\la_\SLC(\Si)$  is defined  (roughly) as the sum of isolated points of
irreducible characters in the $\SLC$ character variety $X(\Si).$
Thus, nontriviality of $\la_\SLC(\Si)$
guarantees the existence of an irreducible representation $\rho\colon \pi_1 \Si \to \SLC$,
and this gives motivation for studying the $\SLC$ Casson invariant.

In this paper, we use the spliced sum construction to present
 a family of  homology spheres with
$\la_\SLC(\Si)=0$.
Since every isolated irreducible character contributes positively to
the $\SLC$ invariant,  homology spheres $\Si$ with $\la_\SLC(\Si)=0$ appear to
be comparatively rare. We prove that, for any
 homology sphere $\Si$ obtained by spliced sum along
two knots in $S^3$, every irreducible representation $\rho\colon
\pi_1 \Si \to \SLC$ lies on a component $X_i$ of the $\SLC$
character variety $X(\Si)$ with $\dim X_i >0$, and this implies
$\la_\SLC(\Si)=0.$

More generally, we investigate the behavior of the invariant $\la_\SLC$
under spliced sum along knots in arbitrary  homology spheres.
Using Casson's surgery formula, Fukuhara and Maruyama, and independently Boyer and Nicas, proved that the $SU(2)$ Casson invariant is additive under spliced sum \cite{FM, BN}.
Unfortunately the same is not always true
for the $\SLC$ Casson invariant.
Counterexamples are provided by Seifert fibered homology spheres.
Recall that $\Si(p,q,r,s)$ is
the spliced sum of $\Si(p,q,rs)$ and $\Si(pq,r,s)$ along the
$rs$-singular fiber in the first and the $pq$-singular fiber in the
second. However,
 $$\la_\SLC(\Si(p,q,r,s)) \neq \la_\SLC(\Si(p,q,rs)) + \la_\SLC(\Si(pq,r,s)).$$
For example, Theorem 2.7 of
\cite{BC} shows that $\la_\SLC(\Si(2,3,5,7))=20$,  whereas
 $\la_\SLC(\Si(2,3,35)) +
  \la_\SLC(\Si(6,5, 7)) = 17+30=47.$

In Theorem \ref{spliced}, our main result,
we develop sufficient conditions, phrased in terms of the knots, under which the
Casson $\SLC$ invariant is additive under spliced sum.

For the remainder of the paper we will use the following notation:
Given a finitely generated group $\pi$, denote by $R(\pi)$ the space
of representations $\rho\colon \pi \to \SLC$ and by $R^*(\pi)$ the
subspace of irreducible representations.  The character of a
representation $\rho$  will be denoted by $\chi_\rho$. The variety
of characters of $\SLC$ representations is denoted $X(\pi)$. Recall
that there is a canonical projection $ R(\pi) \to X(\pi)$ defined by
$ \rho \mapsto \chi_\rho$ which is surjective. Let $X^*(\pi)$ be the
subspace of characters of irreducible characters. Given a manifold
$\Si$, we denote by $R(\Si)$ the space of $\SLC$ representations of
$\pi_1 \Si$ and by $X(\Si)$ the associated character variety.

For the definition of $\la_\SLC$, see \cite{C}.


In section 2 we study  homology spheres resulting from $1/q$ Dehn
surgery on small knots in $S^3$ and show that the $\SLC$ Casson
invariants of such  homology spheres are almost always nontrivial.
In section 3, we introduce splicing and describe the behavior of the
$\SLC$ Casson invariant under spliced sum.

We thank the referee for suggestions improving Theorem 2.1.

\section{Nonvanishing theorems}
In this section, we show that the $\SLC$ Casson invariant is nonzero
for many  homology spheres. Given a knot $K$ in $S^3$ and slope $p/q
\in \QQ \cup \{1/0\}$, we denote by $S^3_{p/q}(K)$ the 3-manifold
obtained by performing $p/q$ Dehn surgery along $K$. Recall that
$S^3_{1/q}(K)$ is always a  homology sphere.

\begin{theorem}
Let $K$ be a small nontrivial knot in $S^3$, and let $q$ be an
integer with $|q|>1$. Then $\la_\SLC(S^3_{1/q}(K)) > 0$.
\end{theorem}

\begin{proof}
By \cite{KM1}, there is an irreducible $SU(2)$ representation of
$\pi_1(S^3_{1/q}(K))$, so the variety of characters of irreducible
$\SLC$ representations is nonempty.  We must show that it contains a
component of dimension 0. In fact we show every component has
dimension 0.

Suppose $q$ is an integer such that the character variety
$X(S^3_{1/q}(K))$ contains a component $Y$ of dimension at least 1.
We may view  $Y$ as a subset of the character variety $X(N)$ of the
complement $N$ of $K$ in $S^3$ since $X(S^3_{1/q}(K)) \subset X(N)$.
Since $K$ is small, $Y$ is one-dimensional, and there is a
well-defined Culler-Shalen seminorm $\|\cdot\|_Y$ on $Y$ given by
\[\|\alpha\|_Y = deg(\tilde{I}^Y_{e(\alpha)}-2)\]
where $\tilde{Y}$ is a smooth projective curve birationally
equivalent to $Y$, $\tilde{I}^Y_{\gamma}$ is the function on
$\tilde{Y}$ induced by the regular function $Y\rightarrow \Bbb C$
taking a character $\xi$ to $\xi(\gamma)$, and $e:H_1(\partial
N;\Bbb Z)\rightarrow \pi_1(\partial N)$ is the inverse of the
Hurewicz isomorphism.  But $\tilde{I}^Y_{e(1/q)} - 2$ vanishes on
$Y$ since $Y$ lies in $X(S^3_{1/q}(K))$, so $Y$ is a $r$-curve as
defined in \cite{BZ1} with $r=1/q$.  Then by Corollary 6.7 of
\cite{BZ1}, we see that $1/q$ is an integer, contradicting the
assumption that $|q|>1.$

It follows that every component of $X(S^3_{1/q}(K))$ has dimension
0, whence the theorem.
\end{proof}

In particular, we have the following:

\begin{theorem} \label{pos2}
If $K$ is a 2-bridge knot or a torus knot, then
$\la_\SLC(S^3_{1/q}(K)) > 0$ for all nonzero integers $q$.
\end{theorem}

\begin{proof} If $|q|>1$, the claim follows from the previous
theorem. By \cite{KM1} we know that $X(S^3_{\pm 1}(K))$ contains an
irreducible character. We show that every component of $X(S^3_{\pm
1}(K))$ is 0-dimensional, so $X(S^3_{\pm 1}(K))$ contains an
isolated irreducible character.

Now as above we know that if $Y$ is a component of dimension greater
than 1 in $X(S^3_{\pm 1}(K))$, then $Y$ has dimension 1, and the
Culler-Shalen seminorm associated to $Y$ is indefinite with $\|\pm
1\|=0$.

It follows by Proposition 5.4 of \cite{BZ3} that there is a positive
integer $k$ and an integral boundary slope $\alpha$ for $K$ such
that the Culler-Shalen seminorm for the curve $Y$ is given by
\[4 \|p \mer + q \lng \|_Y = k|p - q \alpha|\] for any slope
$p/q$. If $\| 1\|=0$, we see that $\alpha = 1$, and if $\|-1\|=0$,
then $\alpha = -1$. But the boundary slopes of 2-bridge knots are
all even integers, and the boundary slopes of the $(r,s)$-torus knot
are 0 and $rs$. Thus in neither case is either 1 or -1 a boundary
slope, so no such curve $Y$ exists.

Thus, $X(S^3_{\pm 1}(K))$ contains an irreducible character and
contains only 0-dimensional components. Hence
$\la_\SLC(S^3_{1/q}(K)) > 0$.
\end{proof}

\section{Splicing}
 The goal of this section is to investigate the behavior of the $\SLC$ Casson invariant
under the operation of spliced sum. Suppose $K_1$ and $K_2$ are
knots in closed 3-manifolds $\Si_1$ and $\Si_2,$ respectively, and
let $M_1 = \Si_1 \setminus K_1$ and $M_2 = \Si_2 \setminus K_2$
denote their complements. Both $M_1$ and $M_2$ are manifolds with
boundary $\partial M_1 = \partial M_2 =T$ a torus, and we denote by
$\mer_i$ and $ \lng_i$ the meridian and longitude of $K_i$ for
$i=1,2.$ The {\it spliced sum} of $K_1$ and $K_2$  is the 3-manifold
$\Si = M_1 \cup_{T} M_2$, with $\partial M_1$ glued to $\partial
M_2$ by a diffeomorphism identifying $\mer_1$ to $\lng_2$ and
$\lng_1$ to $\mer_2.$ If $\Si_1$ and $\Si_2$ are both  homology
spheres, then an elementary exercise shows that $\Si$ is also a
homology sphere.

\begin{figure}[h]
\begin{center}
\psfrag{X}{$M_1$} \psfrag{Y}{$M_2$} \leavevmode\hbox{}
\includegraphics[height=1.0in]{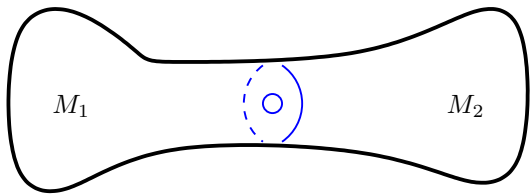}
\caption{{The spliced sum $\Sigma$ along two knots $K_1$ and
$K_2$.}} \label{splice}
\end{center}
\end{figure}

Given a representation $\rho \colon \pi_1 \Si \to \SLC$, then by restriction
we obtain representations $\rho_1 =\rho|_{\pi_1(M_1)}$ and $\rho_2=\rho|_{\pi_1(M_2)}$.
The next theorem shows that any irreducible character $\chi_{\rho}$ for which both
of the induced characters $\chi_{\rho_1}$ and $\chi_{\rho_2}$ are
irreducible must lie on a curve of characters.  Therefore such
representations do not contribute to $\la_{\SLC}(\Si)$.

\begin{theorem} \label{curve}
Suppose $\Si$ is the spliced sum of two 3-manifolds $\Si_1$
and $\Si_2$
and $\chi_\rho \in X(\Si)$ is the character of an irreducible
representation $\rho \colon \pi_1 \Si \to \SLC$ for which the
induced characters $\chi_{\rho_1}$ and $\chi_{\rho_2}$ are
irreducible. Then $\chi_\rho$ lies on a component $X_i$ of $X(\Si)$ with
$\dim X_i
>0.$
\end{theorem}

\begin{proof}
By the Seifert-Van Kampen Theorem, any two irreducible representations
$\rho_1\colon \pi_1 (M_1) \to \SLC$ and $\rho_2\colon \pi_1 (M_2)
\to \SLC$ determine an irreducible representation $\rho\colon \pi_1
\Si \to \SLC$ provided the splicing relations holds, i.e.~provided
that $\rho_1(\mer_1) = \rho_2(\lng_2)$ and  $\rho_2(\mer_2) =
\rho_1(\lng_1).$  We use this fact to construct a curve of
characters in the character variety containing $\chi_{\rho}$.

Since $\rho_1$ and $\rho_2$ are irreducible, they both have
stabilizer subgroup under the conjugation action the group of
central matrices $\left\{
\begin{pmatrix} \pm 1 & 0 \\ 0 & \pm 1
\end{pmatrix} \right\}.$ On the other hand, the restriction $\rho
|_{\pi_1 T}$ of $\rho$ to the splice torus is abelian. Hence its
stabilizer subgroup $\Gamma = {\rm Stab}(\rho|_{\pi_1 T})$ is either
the subgroup $\left\{ \left. \begin{pmatrix} a & 0 \\ 0 & a^{-1}
\end{pmatrix} \right| a \in \CC^* \right\}$ of diagonal matrices or
the subgroup
 $\left\{ \left. \begin{pmatrix} \pm 1 & a \\ 0 & \pm 1 \end{pmatrix} \right| a \in \CC \right\}$
 of upper triangular univalent matrices. In either case, $\dim \Gamma = 1.$
 For any element $\ga \in \Gamma$, the pair $(\rho_1, \ga \rho_2 \ga^{-1})$ is a pair
 of irreducible representations of $\pi_1 (M_1)$ and $\pi_1 (M_2)$ that satisfy
 the splicing relations. The association $\ga \in \Gamma \to \rho_\ga$  gives
 a one-parameter family $\rho_\ga$ of $\SLC$ representations of $\pi_1 \Si$, and it is
 not hard to check that $\rho_\ga$ is conjugate to $\rho$ if and only if $\ga = \pm I.$
 Since distinct conjugacy classes of irreducible representations determine
 distinct characters, this shows that $\chi_\rho$ lies on a component $X_i$ of
 irreducible characters with $\dim X_i >0.$
\end{proof}

The next several results rely on Proposition 6.1 of \cite{CCGLS}
regarding the complement $M_1$ of a knot $K_1$ in a homology sphere
$\Si_1$ . This result asserts that the fundamental group of $M_1$
has a nonabelian reducible representation into $\SLC$ with
eigenvalue $\mu$ if and only if $\mu^2$ is a root of the Alexander
polynomial. Note that in this case the representation has the same
character as an abelian representation, so such characters are the
points of intersection of the curve of reducible characters with
$X^*(\Si_1)$, as is noted in Proposition 6.2 of the same paper. For
any knot $K$ in a homology sphere, let $\Delta_{K_i}(t)$ denote the
Alexander polynomial.

 \begin{proposition}  \label{ind-irr}
Given knots $K_1 \subset \Si_1$ and $K_2 \subset \Si_2$  in homology
spheres, denote their complements $M_1 = \Si_1 \setminus K_1$ and
$M_2 = \Si_2 \setminus K_2$. If $\rho\colon \pi_1 \Si \to \SLC$ is
an irreducible representation of the spliced sum $\Si = M_1
\cup_{T^2} M_2$, then at least one of $\rho_1 =\rho|_{\pi_1 M_1}$ or
$\rho_2= \rho|_{\pi_1 M_2}$ is irreducible.
\end{proposition}
\begin{proof}
We will prove that if $\rho_1$ and $\rho_2$ are both reducible, then
$\rho$ is trivial. Since  $\lng_1$ lies in the second derived
subgroup of $\pi_1(M_1),$ reducibility of $\rho_1$ gives that
$\rho_1(\lng_1) = I.$ Similarly, if $\rho_2$ is reducible, then
$\rho_2(\lng_2) = I.$ Combined with the splicing relations, these
facts imply that $\rho_1(\mer_1) = I = \rho_2(\mer_2)$. Now
$\Delta_{K_1}(1)=\pm 1 \neq 0$, so Proposition 6.1 of \cite{CCGLS}
shows that $\rho_1$ is abelian. Since $\rho_1$ is abelian, it
factors through $H_1(M_1)$, and hence $\rho_1(\mer_1)=I.$ This
implies $\rho_1$ is trivial. A similar argument shows that $\rho_2$
is trivial; hence $\rho$ is trivial.
\end{proof}

A direct consequence is that, for spliced sums along two knots $K_1$ and $K_2$
in  $S^3$, the
$\SLC$ Casson invariant vanishes.

\begin{corollary} \label{vanishing}
If $\Si$ is a spliced sum along two knots $K_1$ and $K_2$ in
$S^3$, then $\la_{\SLC}(\Si) = 0$.
\end{corollary}

\begin{proof}
Suppose $\rho$ is an irreducible representation of $\pi_1(\Si)$ in
$\SLC$ with restrictions $\rho_1$ and $\rho_2$ as before. If
$\rho_1$ and $\rho_2$ are both irreducible, then Theorem \ref{curve}
shows that $\chi_\rho$ is not isolated and hence does not contribute
to $\la_\SLC(\Si).$ Otherwise, by Proposition \ref{ind-irr}, exactly
one of  $\rho_1$ and $\rho_2$ is irreducible.

Suppose that $\rho_1$ is irreducible and $\rho_2$ is reducible. The
reducibility of $\rho_2$ implies $\rho_2(\lng_2)=I,$ so
$\rho_1(\mer_1) = I$ by the splicing relation. However, since $K_1$
is a knot in $S^3,$ we know that the meridian $\mer_1$ normally
generates $\pi_1 (M_1)$. Thus $\rho_1(\mer_1) = I$ implies that
$\rho_1$ is trivial, contradicting the irreducibility of $\rho_1$. A
similar argument with the roles of $\rho_1$ and $\rho_2$ reversed
reveals that $X^*(\Si)$ does not contain any components of dimension
zero.  Therefore $\la_{\SLC}(\Si)=0$.
\end{proof}

The next theorem is our main result, asserting additivity of the
$\SLC$ Casson invariant for spliced sums under certain restrictions.
The restrictions we impose are necessary to rule out the types of
counterexamples that were presented in the introduction.
Specifically, the conditions given below use Proposition 6.1 of
\cite{CCGLS} to rule out unwanted interplay between the reducible
and irreducible characters of $M_1$ and $M_2$.

Before stating the theorem, we find it convenient to define
$$X^\bullet(\Si) = \{ \chi_\rho \in X^*(\Si) \mid \chi_\rho \text{ is isolated}\}$$
to be the subset of {\sl isolated} irreducible characters of $\pi_1(\Si)$.

\begin{theorem} \label{spliced}
Assume $K_1 \subset \Si_1$ and $K_2 \subset \Si_2$ are knots in homology
spheres, and consider the following conditions:
\begin{enumerate}
\item[(i)] For  $\chi_{\rho} \in X^*(\Si_1)$,
if $\mu$ is an eigenvalue of $\rho(\lng_1)$,
then $\Delta_{K_2}(\mu^2)\neq 0$.
\item[(ii)] 
For $\chi_{\rho} \in X^*(\Si_2)$,
if $\mu$ is an eigenvalue of $\rho(\lng_2)$,
then $\Delta_{K_1}(\mu^2)\neq 0$.
\end{enumerate}
If condition (i) is satisfied for all $\chi \in X^\bullet(\Si_1)$
and condition (ii) is satisfied for all $\chi \in X^\bullet(\Si_2)$,
then for the spliced sum, we have
$$\la_{\SLC}(\Si) =  \la_{\SLC}(\Si_1)+ \la_{\SLC}(\Si_2).$$
\end{theorem}

\begin{proof}
If $\rho \colon \pi_1 \Si \to \SLC$ is an irreducible
representation, then Proposition \ref{ind-irr} implies that one of
$\rho_1$ or $\rho_2$ is irreducible. If in addition $\chi_{\rho} \in
X^\bullet(\Si)$ is isolated, then Theorem \ref{curve} shows that
exactly one of $\rho_1$ and $\rho_2$ is irreducible. Hence we can
partition $X^\bullet(\Si)=X^\bullet_1 \cup X^\bullet_2$, where
$$X^\bullet_1 = \{ \chi_\rho \mid  \text{$\rho_1$ is irreducible and  $\rho_2$ is reducible} \},$$
$X^\bullet_2$ is defined similarly, and $X^\bullet_1$ and $X^\bullet_2$ are disjoint.

For $\chi_\rho \in X^\bullet_1,$ reducibility of $\rho_2$ and the
splicing relations imply that $\rho_1(\mer_1) = \rho_2(\lng_2) = I.$
Hence $\rho_1$ extends to an irreducible representation
$\rho_1'\colon \pi_1(\Si_1) \to \SLC.$ Thus, we have a natural map
$\Phi_1: X^\bullet_1 \to X^*(\Si_1)$ given by $\chi_\rho \mapsto
\chi_{\rho_1'}$. We define a map $\Phi_2: X^\bullet_2 \to
X^*(\Si_2)$ analogously.

Conversely, given  an irreducible representation $\rho_1' \colon
\pi_1(\Si_1) \to \SLC$ with $\chi_{\rho'_1}$ isolated,  we define a
reducible representation $\rho_2 \colon \pi_1(M_2) \to \SLC$ by
setting $\rho_2(\mer_2) = \rho_1(\lng_1)$. (Here, $\rho_1 =
\rho_1'|_{\pi_1(M_1)}.$) Note that hypothesis (i) implies that
$\rho_2$ is abelian by Proposition 6.1 of \cite{CCGLS}, so this
assignment of $\rho_2(\mer_2)$ completely determines $\rho_2$.
Direct inspection shows that $\rho_1$ and $\rho_2$ satisfy the
splicing relations; thus they give rise to an irreducible
representation $\rho \colon \Si \to \SLC$. Further, since
$\chi_{\rho_1'}\in X^\bullet(\Si_1)$ is isolated and $\rho_2$ is
completely determined by $\rho_1(\lng_1),$ it is not difficult to
see that $\chi_\rho \in X^\bullet_1 \subset X^\bullet(\Si)$ is also
isolated.

This defines a map $\Psi_1 \colon X^\bullet(\Si_1) \to X^\bullet_1,$
which is an inverse to $\Phi_1$ and gives a one-to-one
correspondence between $X^\bullet(\Si_1)$ and $X^\bullet_1$.

The same construction with the roles
of $\rho_1$ and $\rho_2$ reversed defines
a map $\Psi_2 \colon X^\bullet(\Si_2) \to X^\bullet_2$ which is an
inverse to $\Phi_2$ and gives a one-to-one correspondence between
$X^\bullet(\Si_2)$ and $X^\bullet_2$.

It remains to show that $\chi_1 \in X^\bullet_1\subset
X^\bullet(\Si)$ and $\Phi_1(\chi_1) \in X^\bullet(\Si_1)$ both
contribute equally to their respective $\SLC$ Casson invariants, and
similarly for $\chi_2 \in X^\bullet_2\subset X^\bullet(\Si)$ and
$\Phi_2(\chi_2) \in X^\bullet(\Si_2).$

Choose a triangulation of $\Si_1$ such that the 1-skeleton contains
$K_1$. Build a Heegaard decompostion $(U_1,U_2)$ of $\Si_1$ by
letting $U_1$ be a tubular neighborhood of this 1-skeleton. (See
Theorem 2.5 of \cite{H}.) Call the Heegaard surface for this
Heegaard decomposition $F_1$.  Similarly choose a triangulation for
$\Si_2$ whose 1-skeleton contains $K_2$ and build a Heegaard
splitting $(V_1,V_2)$ of $\Si_2$ by letting $V_2$ be a neighborhood
of the 1-skeleton.  Call the Heegaard surface $F_2$.

Choose a symplectic basis for $F_1$ consisting of curves $a_1,
b_1, \ldots, a_j, b_j, \mer_1, \lng_1,$ where the curves
$a_1,\ldots, a_j,$ and $\mer_1$ are homotopically trivial in
$U_1.$ Choose a symplectic basis for $F_2$ consisting of curves
$\lng_2, \mer_2, c_1, d_1, \ldots, c_k, d_k,$  where the
curves $\mer_2$ and $d_1, \ldots, d_k$ are homotopically trivial
in $V_2.$

Note that $U_1$ is the union of a tubular neighborhood $N(K_1)$ of
the knot $K_1$ and a handlebody $H_1$ of genus $j$ spanned by  the
curves $a_1, b_1, \ldots, a_j, b_j$ in the obvious way. Similarly
$V_2$ is the union of a tubular neighborhood $N(K_2)$ of the knot
$K_2$ and a handlebody $H_2$ of genus $k$ spanned by the curves
$d_1, c_1, \ldots, d_k, c_k$ in the obvious way. On the other hand,
$U_2$ is a subset of $M_1$ and $V_1$ is a subset of $M_2$. From this
one sees that the restrictions of these Heegaard splittings to $M_1$
and $M_2$ glue together to form a Heegaard splitting $(W_1,W_2)$ of
$\Si$, where $W_1$ is the connected sum of $H_1$ and $V_1$ and $W_2$
is the connected sum of $U_2$ and $H_2$. Denote by $F$ the Heegaard
surface of this Heegaard decomposition of $\Si.$  Note that $F$ can
be viewed as the connected sum of $F_1$ and the boundary of $H_2$ or
as the connected sum of the boundary of $H_1$ and $F_2$.

Now given $\chi_{\rho}\in X^\bullet_1$, we see that $X(V_2 -
N(K_2))$ and $X(V_1)$ are transverse in $X(F_2)$ at $\chi_{\rho_2}$
since the dimension of $H^1(M_2; \slc_{\Ad \rho})$ is 1. This follows from
the Mayer-Vietoris sequence for $M_2 = (V_2 - N(K_2)) \cup V_1$,
using the fact that condition (i) of the theorem is satisfied at
$\rho$.

It follows that there is an isotopy $h_t$ of $X(F)$, $t\in [0,1]$,
such that $h_0$ is the identity; $h_t(\phi(c_i)) = \phi(c_i)$ ,
$h_t(\phi(d_i)) = \phi(d_i)$, $h_t(\chi_{\phi}(\mer_2)) =
\chi_{\phi}(\mer_2)$, and
$h_t(\chi_{\phi}(\lng_2))=\chi_{\phi}(\lng_2)$ for every
$\chi_{\phi}$ and every $i$; and $h_1(X(W_2))$ meets $X(W_1)$
transversely in a neighborhood of $\chi_{\rho}$.  In fact, $h_t$ can
be chosen to have support in a neighborhood $N$ of $\chi_{\rho}$
such that $N$ meets $X(\Si)$ only in $\chi_{\rho}$ and for any
$\chi_{\phi}$ in $N$, if $\mu$ is an eigenvalue of $\phi(\lng_1)$,
then $\Delta_{K_2}(\mu^2)\neq 0$. Then the contribution of
$\chi_{\rho}$ to $\la_{\SLC}(\Si)$ is precisely the number of points
in the intersection of $X(W_1)$ and $h_1(X(W_2))$ in $N$.

 Now since $h_t(\phi(c_i)) = \phi(c_i)$ and
$h_t(\phi(d_i)) = \phi(d_i)$ for every i and every $\phi$, we see
that $h_t$ preserves the subvariety of $X(F)$ consisting of
characters $\chi_{\phi}$ for which
$[\phi(c_1),\phi(d_1)][\phi(c_2),\phi(d_2)]\ldots[\phi(c_k),\phi(d_k)]=I$.
But these characters are precisely the characters for which
$[\phi(a_1),\phi(b_1)][\phi(a_2),\phi(b_2)]\ldots[\phi(a_j),\phi(b_j)][\phi(\mer_1),\phi(\lng_1)]=I$
- i.e. the characters which are the images of characters in
$X(F_1)$. It follows that $h_t$ induces an isotopy $\tilde{h}_t$ of
$X(F_1)$. Moreover $X(U_1)$ and $\tilde{h}_1(X(U_2))$ intersect
transversely in a neighborhood of the image of $\Phi_1(\chi_{\rho})$
in $X(F_1)$ since $X(W_1)$ and $h_t(X(W_2))$ are transverse.  Thus
the contribution of $\Phi_1(\chi_{\rho})$ to $\la_{\SLC}(\Si_1)$ is
precisely the number of points of intersection of $X(U_1)$ and
$\tilde{h}_1(X(U_2))$ in the neighborhood of $\Phi_1(\chi_{\rho})$
which is the support of $\tilde{h}$.

It remains to be shown that the points of intersection of $X(U_1)$
and $\tilde{h}_1(X(U_2))$ in the neighborhood of
$\Phi_1(\chi_{\rho})$ which is the support of $\tilde{h}$ are in
one-to-one correspondence with the points of intersection of
$X(W_1)$ and $h_1(X(W_2))$ in $N$. This follows since every point
$\chi_{\phi}$ in the intersection of $X(W_1)$ and $h_1(X(W_2))$ in
$N$ satisfies
 $\phi(d_1) = \phi(d_2)=\ldots\phi(d_k)=I$ since $h_t$
did not affect the values of $\phi$ at $d_1,d_2,\ldots,d_k$ and
$\chi_{\phi}\in h_t(X(W_2))$, and so
$[\phi(c_1),\phi(d_1)][\phi(c_2),\phi(d_2)]\ldots[\phi(c_k),\phi(d_k)]=I
=
[\phi(a_1),\phi(b_1)][\phi(a_2),\phi(b_2)]\ldots[\phi(a_j),\phi(b_j)][\phi(\mer_1),\phi(\lng_1)]$.

Thus, $\chi_{\rho}$ and $\Phi_1(\chi_{\rho})$ contribute equally to
their respective Casson invariants. That points in $\chi_2 \in
X^\bullet(\Si_2)$ and $\Psi_2(\chi_2) \in X^\bullet_2 \subset
X^\bullet(\Si)$ also contribute equally to their respective Casson
invariants can be proved analogously.
\end{proof}

\begin{remark}
A useful observation is that the two hypotheses in Theorem \ref {spliced},
are equivalent to the following conditions:
\begin{enumerate}
\item[(i)] If $t^{2k}$ is a root of the Alexander polynomial of $K_2$,
then $A_{K_1}(t,t^{-k}) \neq 0$, where $A_{K_1}$ denotes the
$A$-polynomial of $K_1$ .
\item[(ii)]
If $t^{2k}$ is a root of the Alexander polynomial of $K_1$, then
$A_{K_2}(t,t^{-k}) \neq 0$, where $A_{K_2}$ denotes the $A$-polynomial
of $K_2$ .
\end{enumerate}

\end{remark}

We now describe the operation of $k$-spliced sum for
two knots $K_1, K_2$ in $S^3$. Let
$M_1 = S^3 \setminus K_1$ and $M_2 = S^3 \setminus K_2$ be their
complements, and denote by $\mer_1, \lng_1$ and $\mer_2, \lng_2$ the
meridian and longitude of $K_1$ and $K_2$.
The {\it $k$-spliced sum} is the 3-manifold $\Si_k = M_1 \cup_{\phi}
M_2$, with $\partial M_1$ glued
to $\partial M_2$ by a diffeomorphism $\phi$ identifying $\mer_1$ to
$\lng_2$ and $\lng_1$ to $\mer_2 \lng_2^k.$
The diffeomorphism $\phi \colon \partial M_1 \to \partial M_2$
is represented on  $\pi_1 T$ by the matrix
$\begin{pmatrix} 0 & 1 \\ 1 & k \end{pmatrix}.$
 It is not difficult to see that
$\Si_k$ is an  homology sphere with $\Si_0$ the spliced
sum considered previously. Further, if $K_1$ is the unknot, then
$\Si_k$ is the  homology sphere obtained by $1/k$ Dehn surgery on
$K_2.$ We apply Theorem \ref{spliced} to determine the $\SLC$ Casson invariant
of $k$-spliced sums.

\begin{corollary}
Let $K_1$ and $K_2$ be knots in $S^3$.  Let $k$ be an integer, and
let $\Si_k$ be the $k$-spliced sum of $K_1$ and $K_2$. Then
$\la_{\SLC}(\Si_k) = \la_\SLC  \! \left(S^3_{1/k} \left(K_2\right)\right)$.
\end{corollary}

\begin{proof} Set $\Si_2 =S^3_{1/k}(K_2)$ and  let $\widetilde{K}_2$ be
the image of $K_2$ in $\Si_2$. Then the $k$-spliced sum of $K_1$ and $K_2$ is
homeomorphic to the spliced sum of $S^3$ and $\Si_2$ along $K_1$ and $\widetilde{K}_2$
since the meridian of $\widetilde{K}_2$
is $\mer_2 \lng_2^k$.

If $\rho_2 \colon \pi_1(\Si_2) \to \SLC$ is an irreducible representation,
then $\rho_2(\mer_2 \lng_2^k) = I$.
We can conjugate so that $\rho_2(\lng_2)$ is either diagonal or a matrix of the
form $ \begin{pmatrix} \pm 1 & a \\ 0 & \pm 1
\end{pmatrix}$ for some $a \in \CC$.
If $\rho_2(\lng_2)$
is diagonal, then since $\rho(\mer_2 \lng_2^k) = I$,
we see that $\rho_2(\mer_2)$
is also diagonal, which contradicts
the irreducibility of $\rho_2$.

Hence the eigenvalues of $\rho_2(\lng_2)$ are in $\{\pm
1\}$. Since their squares, which equal 1, are not
roots of the Alexander polynomial for any knot,  Theorem \ref{spliced}
applies and implies $\la_{\SLC}(\Si) = \la_\SLC(\Si_2)=\la_\SLC(S^3_{1/k}(K_2))$.
\end{proof}

Finally, combining this corollary with Theorem \ref{pos2} yields
the following result.

\begin{corollary}
If $K_1$ is any knot in $S^3$ and $K_2$ is a 2-bridge or a torus
knot, then the $k$-spliced sum of $K_1$ and $K_2$  satisfies
$\la_\SLC(\Si_k)
> 0$ for $k \neq 0.$
\end{corollary}

\begin{figure}[h]
\begin{center}
\psfrag{A}{$k=0$} \psfrag{B}{$k=-5$} \leavevmode\hbox{}
\includegraphics[height=1.5in]{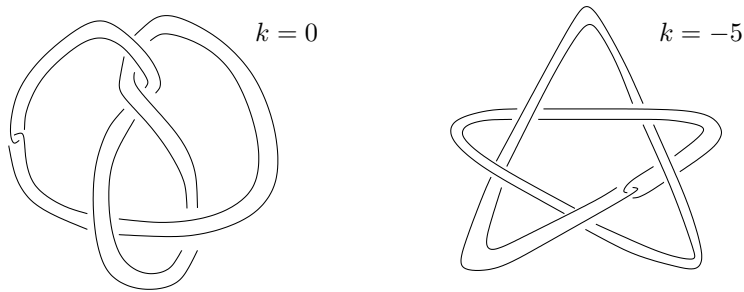}
\caption{{The $k$-twisted Whitehead doubles of the figure-8 knot and
the $(2,5)$ torus knot.}} \label{whitehead}
\end{center}
\end{figure}

If $K_1$ is the left-hand trefoil, the $k$-spliced sum of $K_1$ and
$K_2$ is the  homology sphere obtained by $-1$ surgery on the
$-k$-twisted Whitehead double of $K_2$ (see Prop. 6.1, \cite{KKR}).
In particular, we conclude that the  homology sphere $\Si_k$ obtained
by $-1$ surgery on a $k$-twisted Whitehead double of any 2-bridge or
torus knot has $\la_\SLC (\Si_k)>0$ provided $k \neq 0$. For example,
taking $K_2=T(p,q)$ the $(p,q)$ torus knot and denoting by $L_{k}$
the $-k$-twisted Whitehead double of $T(p,q)$, we see that for $k>0$, we have
$$\la_\SLC( S^3_{-1}(L_{k}))
=  \la_\SLC(\Si(p,q,pqk-1) = \frac{(p-1)(q-1)(pqk-2)}{4}$$ by
combining the above corollary with Theorem 2.3, \cite{BC}. A similar result follows for $k<0$,
and the same
idea applies to provide explicit computations of
 $\la_\SLC( S^3_{-1}(L_{k}))$ for $L_k$  the $-k$-twisted Whitehead double of a
 twist knot, see Theorems 5.7 and 5.9, \cite{BC}.

\bibliographystyle{amsplain}

\end{document}